\documentclass[11pt]{article}

\usepackage{amssymb,latexsym,amsmath,amsthm,multicol,float,amstext, color}
\usepackage [final]{graphicx}
\numberwithin{equation}{section}

\newcommand{\be}{\begin{eqnarray}}
\newcommand{\ee}{\end{eqnarray}}
\newcommand{\ce}{\begin{eqnarray*}}
\newcommand{\de}{\end{eqnarray*}}
\newtheorem{thm}{Theorem}[section]
\newtheorem{lemma}[thm]{Lemma}
\newtheorem{remark}[thm]{Remark}

\begin{document}

\title {The Proof of the Collatz Conjecture }

\author{Chin-Long Wey}

\maketitle
\begin{abstract}
The \textit{3n+1}, or \textit{Collatz problem}, is one of the hardest math problems, yet still unsolved. The \textit{Collatz conjecture} is to prove or disprove that the Collatz sequences COL(n) always eventually reach the number of 1, for all n$\in$N$^+$ (all positive integers). The Syracuse conjecture is a (2N+1)-version of Collatz conjecture, where (2N+1) is all positive odd integers. 

The \textit{Syracuse} and \textit{Collatz} problems can be conceptually described by a \textit{tree trunk} and \textit{branches}. The trunk is made of the \textit{junctions} that produce the main branches, where  $J_0$=1 is the \textit{root junction}. Each branch consists of \textit{active} and \textit{dead} junctions, where only the active junctions are capable of producing new sub-branches.  
Conceptually assuming the trunk and branches can grow indefinitely and can also absorb nutrients from the root. As the tree grows indefinitely, all N$^+$ (2N+1) are included for the Collatz (Syracuse) sequence. 
This paper develops the \textit{inverse Collatz function} and the \textit{inverse Syracuse functions} to construct the tree trunk and branches starting from the root junction $J_0$=1 and assign the also positive (odd) integers to all junctions. To verify the Collatz (Syracuse) sequences always eventually reach the number of 1,  this paper also develops the \textbf{PathFinding} algorithm. Given n$\in$N$^+$ (2N+1), the algorithm finds a path from n to the root junction $J_0$=1 by the virtual tree structure to prove both Syracuse and Collatz conjectures.
\end{abstract}

\section{Introduction}

Let N=$\{$0,1,2,$\cdots\}$, $N^+$=$\{$1,2,$\cdots\}$, 2N+1=$\{$1,3,5,$\cdots\}$, and $2N^+$=$\{$2,4,6,$\cdots\}$.  
The \textit{3n+1 problem}, or \textit{Collatz problem}, is one of the hardest math problems, yet still unsolved [1]. The \textit{Collatz function} is $n_{i+1}$=Col($n_i$), where 
\begin{equation}
Col(n_i) = 
    \begin{cases}
  3*n_i+1, & if \quad n_i \in 2	N+1   \\
       n_i/2, & if \quad n_i \in 2N^+
    \end{cases}
\end{equation}

The Collatz conjecture is to prove or disprove that the sequence COL(n) always eventually reach the number of 1 [1]. 
COL(n)=$\{$n, Col(n), Col$^2$(n), …, Col$^c$(n)$\}$, and $n_i$=Col$^i$(n). If Col$^c$(n)=1, then the sequence COL(n) converges to 1, denoted by COL(n)$\rightarrow$1. Otherwise, the sequence is not converged to 1, denoted by COL(n)$\nrightarrow$1.
\[
\textbf{\textit{Collatz Conjecture}}:  COL(n)\rightarrow1, \forall n\in N^+.
\]

The \textit{Syracuse conjecture} is the (2N+1)-version Collatz conjecture. The Syracuse sequence SYR(n)=$\{$n, Syr(n), Syr$^2$(n), …, Syr$^s$(n)$\}$, $n_i \in$2N+1, and 
n$_{i+1}$=Syr($n_i$), where
\begin{equation}
Syr(n_i) = (3*n_i+1)/2^r.	
\end{equation}
The \textit{Syracuse conjecture} is to prove or disprove that the sequence SYR(n) always eventually reach the number of 1 [1-3], i.e., SYR(n)$\rightarrow$1, or Syr$^s$(n)=1,
\[
\textbf{\textit{Syracuse Conjecture}}:  SYR(n)\rightarrow1, \forall n\in 2N+1.
\]
For any Syracuse sequence SYR(n)=$\{$n=$J_s$,$J_{s-1}$,…, $J_1$,$J_0\}$, where Syr($J_r$)=$J_{r-1}$. If $J_0$=1, then SYR(n)$\rightarrow$1. 
For n=45, the Syracuse sequence Syr(45)=$\{$45,17,13,5,1$\}$, and the Collatz sequence COL(45)=$\{$\textbf{45},136,68,34,$\textbf{17},52,26,$\textbf{13},40,20,10,$\textbf{5},16,8,4,2,$\textbf{1}$\}$. \textit{The Syracuse conjecture is the (2N+1)-version Collatz conjecture}.

\textbf{Inverse Syracuse Function}
\quad  By (1.2), for any n$\in$2N+1,t$\in$N, m=Syr(n)=(3n+1)/2$^r\in$ 2N+1, n=$I_b$(k,t) [3],
where r=2k+2 if b=1; and r=2k+1 if b=5, k$\in$N. Thus, n=(2$^r$m-1)/3 =$I_b$(k,t). For k=0,

\quad \quad \quad \quad \quad $I_b$(0,t)=(2$^r$m-1)/3, where r=2 if b=1, and r=1 if b=5. 

\noindent Let g=$I_b$(0,t) (mod 6), the function G is defined as

\begin{equation}
G(J_s)=J_{s+1} = 
    \begin{cases}
  I_b(0,t), & if \quad g \neq 3, and   \\
  I_b(1,t), & if \quad g = 3
    \end{cases}
\end{equation}
\noindent where $I_b$(1,t)=4*$I_b$(0,t)+1 [3]. Let $E^0$=$\{J_0$,$J_1$,$J_2$, …$\}$=$\{$1,5,13,17,…$\}$, as shown in Figure 1(b), SYR(17)=$\{$1,5,13,17,…$\}$,  
Syr$^3$($J_3$)=Syr$^3$(17)=1, SYR(17)$\rightarrow$1, and SYR(x)$\rightarrow$1, $\forall$x$\in E^0$. G(m)=n, or G($J_r$)=J$_{r+1}$, r$ \in$N. 
By (1.3), if $J_0$=1, the sequence (11…1) is a \textit{trivial cycle} of the Syracuse sequences. 

•   Syr(G(m))=Syr($I_b$(0,t))=Syr((2$^r$m-1/3))=(3((2$^r$m-1/3)+1)/2$^r$=m; and 

•   G(Syr(n)) = G((3n+1)/2$^r$)=(2$^r$((3n+1)/2$^r$)-1)/3=n;

\noindent Thus, the function G is the \textit{inverse Syracuse function}. 

\textbf{Inverse Collatz Function}
For any Collatz sequence COL(n)=$\{$n=$J_c$,$J_{c-1}$,…, $J_1$,$J_0\}$, where Col($J_r$)=$J_{r-1}$, r=c,c-1,…,2,1.  
For any n$\in N^+$=$\{$2N+1$\}\cup\{2N^+\}$, by (1.1), 
m=Col(n) =n/2, if n$\in 2N^+$, and  n=2m, which may cause the inverse function to grow indefinitely. However, by (1.1), m=Col(n)=3n+1, if n$\in$2N+1, let d=n=(m-1)/3, the growth stops if d$\in$2N+1 and  d$\neq$3 (mod 6). This is simply because that if d$\notin$2N+1,  
d=(m-1)/3 is even and  m=3d+1 is odd and contradicts to m=3n+1  is even.Thus, the function H is defined as
\begin{equation}
H(m)=
    \begin{cases}
  d, & if \quad d \in 2N+1 \quad and \quad d \neq 3, and   \\
  2m, & Otherwise.
    \end{cases}
\end{equation}

For example, m=$h_6$=10, d=(m-1)/3=3, $h_7$=20, $h_8$=40, d=(40-1)/3=13$\neq$3, thus, $h_9$=13.
Let $V^0$=$\{h_0$,$h_1$,$h_q$,…, $h_c$,…$\}$=$\{$1,2,4,8,16,5,10,20,40,13,26,…$\}$, as shown in Figure 1(c). Col$^{10}$($h_{10}$)=Col$^{10}$(26)=1, COL(26)$\rightarrow$1, and COL(x)$\rightarrow$1, $\forall$ x$\in V^0$.
Note that H(m)$\neq$3 (mod 6), $\forall$m$\in$2N+1, e.g., if m=17, H(17)=2*17=34, d=(34-1)/3=11$\in$2N+1 and  d$\neq$3, thus H(34)=11=$I_5$(0,1), and $I_5$(1,1)=45. Similar to the inverse Syracuse function, H(m)=$I_b$(0,t), or $I_b$(0,t) if $I_b$(0,t)=3 (mod 6). 
By (1.4), $J_0$=1,  $J_1$=2, and  $J_2$=4, the sequence (1,2,4,1,2,4,…) is a \textit{trivial cycle} of the Collatz sequences, and 

•   If H(m)$\in$2N+1, Col(H(m))=3H(m)+1=m; If H(m)$\in 2N^+$,	Col(H(m))=(H(m)/2=m;

•   If n$\in$2N+1, H(Col(n))=(Col(n)-1)/3=n; and If n$\in 2N^+$, H(Col(n))=H(n/2)=n.

\noindent Thus, the function H is the \textit{inverse Collatz function}.

The tree trunk is made of the \textit{junctions}, J-link $E^0$=$\{J_s\}$, $J_s \in$2N+1,
that produce the main branches, $D^0$=$\{D_s\}$, and the root junction is $J_0$=1. 
\begin {equation}
SYR(J_s)\rightarrow 1 \quad and \quad SYR(D_s)\rightarrow 1, \quad \forall s \in N.	 
\end {equation}
By (1.5), SYR(x)$\rightarrow$1, $\forall$x$\in D^0 \subset$2N+1, but $D^0\neq$2N+1.

Each of the main branches, $D_s$, consists of \textit{active} and \textit{dead junctions}, as shown in Figure 1(b), where only the \textit{active junctions} are capable of producing new sub-branches. Let $F_0$ and $B_0$ be the collections of all dead (black) and active (grey) junctions of $D^0$, respectively. The white dots are the junctions in $E^0$.
Each of active junction $J_{x0}$ of $D^0$ produces the sub-tree trunk $E^1$, as shown in Figure 1(e), and each junction of $E^1$, produces the sub-branches of $D^1$. 
Further, each of active junction $J_{x0}$ of Dr produces the sub-tree trunk, $E^{r+1}$, as shown in Figure 1(e), and each junction of $E^{r+1}$, produces the sub-branches 
of $D^{r+1}$. 
SYR($J_s$)$\rightarrow$1, $\forall J_s \in E^{r+1}$, and 
SYR($I_b$(p,t))$\rightarrow$1, $\forall I_b$(p,t)$\in D^{r+1}$.

   
\begin{figure}
\centering
\includegraphics[width=1\linewidth]{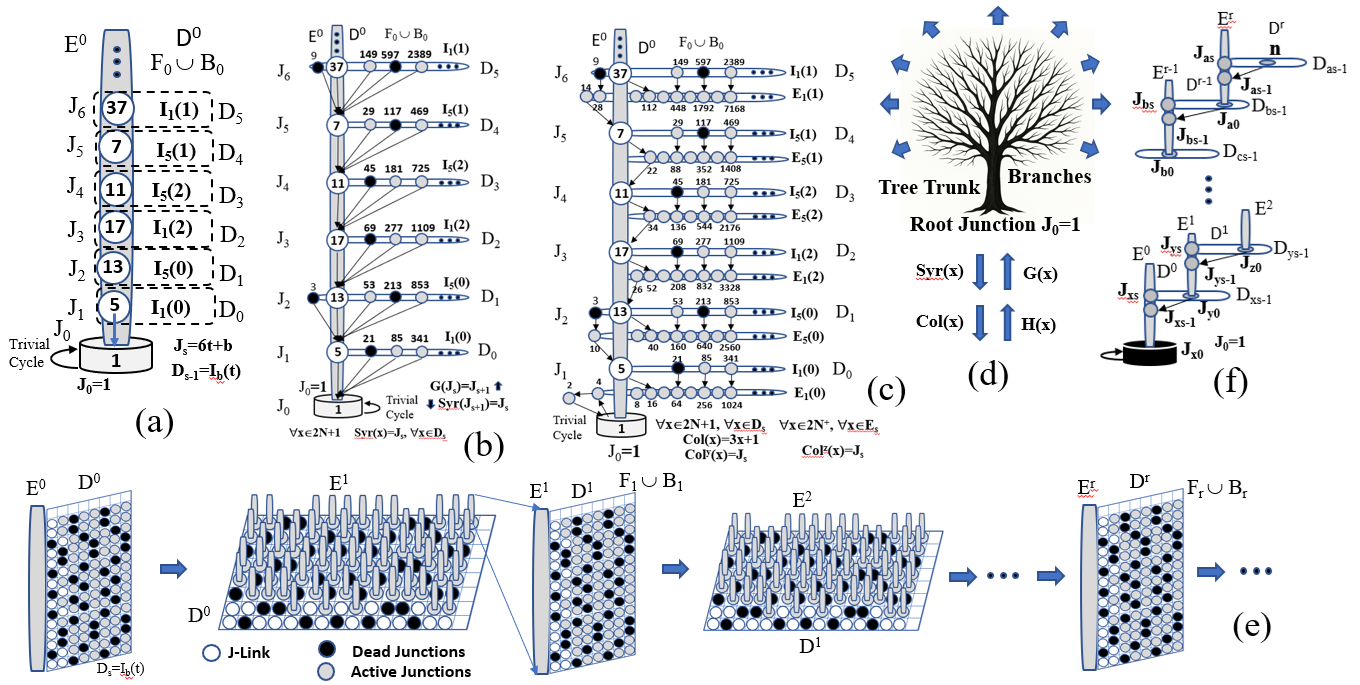}
\caption{\label{fig:IGz}(a) Tree trunk; (b) Tree trunk and Junctions; (c) (2N+1)-version Collatz sequence (or, Syracuse sequence); (e) Expansion of both E$^0$ and $^0$; and (f) Tree trunk structure.}
\end{figure}

\begin {remark} Tree Trunk and Main Branches for Syracuse Sequences
\begin {enumerate}

\item If the tree trunk, sub-tree trunks, and sub-branches can grow indefinitely, then 
SYR($J_s$) $\rightarrow$1, $\forall J_s \in$E, and 
SYR($I_b$(p,t))$\rightarrow$1, $\forall I_b$(p,t)$\in$D, where 
\begin {equation}
\begin {split}
& E=\cup_{r=0}^{\infty}E^r=\{6t+1\}\cup\{6t+5\}; and  \\
& D=\cup_{r=0}^{\infty}D^r=\{I_1(p,q)\}\cup\{I_5(p,q)\}=2N+1.
\end {split}
\end {equation}
and
\begin {equation}
SYR(n) \rightarrow 1, \forall n \in D=2N+1.
\end {equation}

By (1.7), the Syracuse conjecture holds.

\item The tree trunk structure in Figure 1(e) starts from the root junction $J_0$=1, generating $E^0$/$D^0$,  $E^1$/$D^1$, and then generate $E^r$/$D^r$ upward, 
by (1.6), D=2N+1. The junctions $J_s$ of $E^r$ are generated by the \textit{inverse Syracuse function} in (1.3), for all r$\in$N.  For any n$\in$2N+1, n=6t+a$\in D_{as-1}$=$I_a$(t)$\in E^r$. The root junction $J_{a0}$ produces $E^r$($J_{a0}$), one of the sub-tree trunks of $E^r$, and the sub-branch $D_{bs-1}$.. By (1.3), Syr(G($J_{a0}$))=$J_{a0}$.
The same procedure is repeatedly applied until 
$J_{xs}\in E^{r-u}$ and  Syr(G($J_{x0}$))=$J_{x0}$=1, i.e., $E^{r-u}$=$E^0$, or r=u, 
meaning that n locates at $E^u$, and SYR(n)$\rightarrow$1. The Syracuse conjecture holds.
\end {enumerate}
\end {remark}

\section {Properties of Matrices $I_a$(p,q) and $E_b$(p,q)}

The properties of these matrices are summarized as follows, and and shown in Table 1, 

\begin {thm} $\{I_a$(p,q)$\}$, a=1,5, p,q$\in$N [3].

(1) $I_1$(0,q)=8q+1, $I_5$(0,q)=4q+3, and $I_a$(p+1,q)=4*$I_a$(p,q)+1; 

(2) $I_1$(p,q)=[(6q+1)*$4^{p+1}$-1]/3, and $I_5$(p,q)=[(6q+5)*$4^{p+1}$-2]/6. 

(3) $\{I_1$(q)$\}\cup\{I_5$(q)$\}$=2N+1, and the values of all entries are distinct; 

(4) Syr($I_b$(p,q))=Syr($I_b$(0,q))=6q+b.
\end {thm}

\begin {thm} $\{E_b$(p,q)$\}$, b=1,3,5, p,q$\in$N [3].

(1) $E_b$(r,q)=(6q+b)*$2^{r+1}$; r=8q+1 if b=1, or r=4q+3 if b=5;

(2) $\{E_1$(2p,q)$\}\cup\{E_5$(2p+1,q)$\}$=$\{$6q+2$\}$; $\{E_1$(2p+1,q)$\}\cup\{E_5$(2p,q)$\}$=$\{$6q+4$\}$; 

\quad \quad and $\{E_3$(p,q)$\}$=$\{$6q$\}$-$\{$0$\}$;

(3) $\{E_1$(p,q)$\}\cup\{E_3$(p,q)$\}\cup\{E_5$(p,q)$\}$=2$N^+$, and all entries have distinct values; 
\end {thm}

\begin{table}
\centering
\caption{\label{tab:widgets1}Matrices $\{I_1$(p,q)$\}$ and $\{I_5$(p,q)$\}$.}
\includegraphics[width=1.0\linewidth]{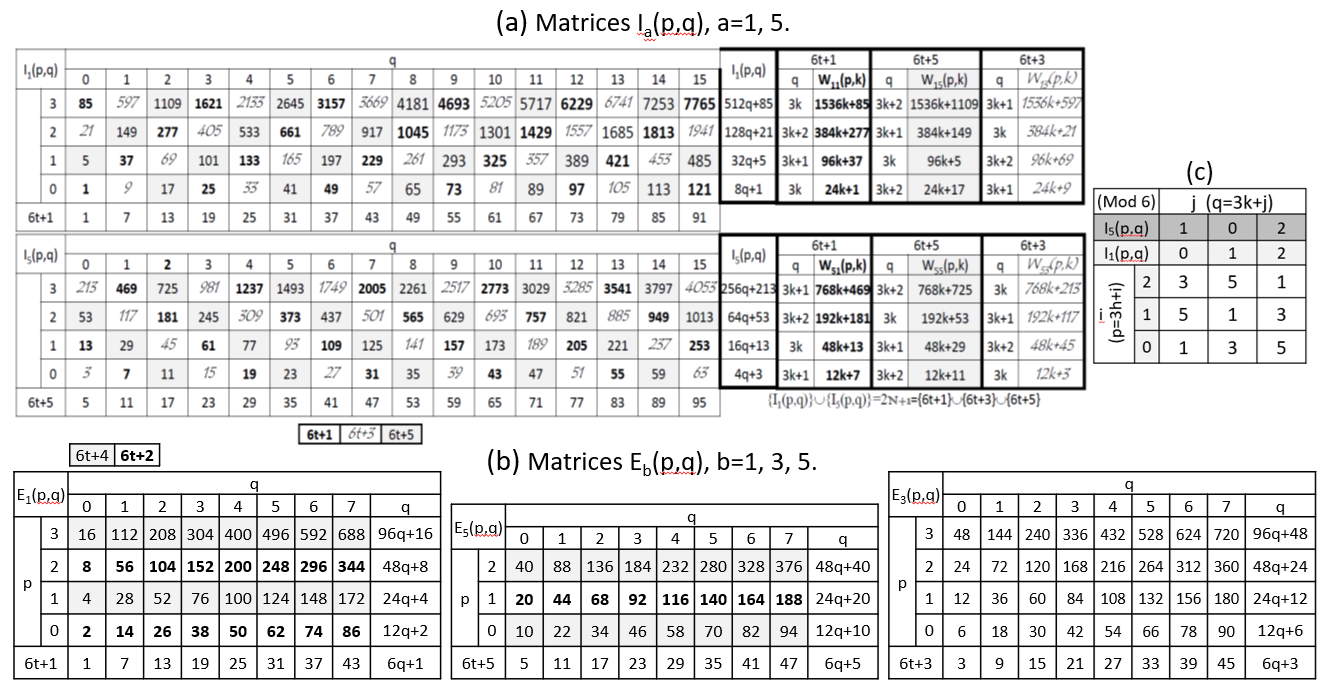}
\end{table}

For the Syracuse conjecture, by Theorem 2.1, each column is defined as 
\begin {equation}
I_a(q)=\{m\in 2N+1 | Syr(m)=6q+a\}	
\end {equation}
In Table 1(a), if a=1, q=2, the column $I_a$(2)=$\{m_0$,$m_1$,$m_2$,$m_3$,…$\}$ =$\{$17,69,277,1109,…$\}$. By Theorem 2.1, each column is defined by (2.1), where $\{I_1$(p,q)$\}\cup\{I_5$(p,q)$\}$=2N+1=
$\{$6t+1$\}\cup\{$6t+3$\}\cup\{$6t+5$\}$.
Let b=$I_a$(p,q) (mod 6), b=1,3,5, the sets $\{$6t+b$\}$ equally share the entries of $\{I_1$(p,q)$\}\cup\{I_5$(p,q)$\}$, marked by boldface, italics, and shade, respectively. The entries with b=1, 3, or 5 are periodically appeared to both row- and column-directions of the matrices, as shown in Table 1(a); 
Let $u^p$=$I_a$(p,q), $v^p$=$u^p$ (mod 6), p=0,1,2. By Table 1(c), if $v_0$=1, then $v_1$=5, and $v_2$=3; if $v_0$=3, then $v_1$=1, and $v_2$=5; and if $v_0$=5, then $v_1$=3, and $v_2$=1. By Theorem 2.1(4), Syr($I_b$(p,q))=Syr($I_b$(0,q))=6q+b.

Let $W_{ab}$(p,q)=$I_a$(p,q), $\{W_{ab}$(p,q)$\}$=$\{I_a$(p,q)$\}\cup\{$6t+b$\}$, and $\{W_{1b}$(p,q)$\}\cup\{W_{1b}$(p,q)$\}$= $\{$6t+b$\}$, b=1,3,5. For any n=6t+b, there exists $I_a$(p,q), such that $I_a$(p,q)=$W_{ab}$(p,q)=6t+b [3], SYR($I_a$(p,q))= SYR(6t+b), and 
\begin {equation}
SYR(I_a(p,q))\rightarrow 1,\quad  iff \quad SYR(6t+b)\rightarrow 1. 
\end {equation}
For any $I_a$(p,q)$\in$2N+1, $\forall I_b$(t)$\in \{I_1$(q)$\}\cup\{I_5$(q)$\}$, and 
\begin {equation}
SYR(I_a(p,q))\rightarrow 1, \quad iff \quad SYR(I_b(t))\rightarrow 1.
\end {equation}

By Theorem 2.1(3), $\{I_1$(p,q)$\}\cup\{I_5$(p,q)$\}$=2N+1, and the values of all entries in Table 1(a) are distinct. 
The columns $I_a$(q) of $I_a$(p,q), a=1,5, are arranged in terms of 6q+a, q$\in$N. By the G-function in (1.3),
the root junction $J_0$=1 produces the J-link, $E^0$ and the main branches, $D^0$, 
where $E^0$=$\{$1,5,13,17,…$\}$,
$D^0$=$\{I_1$(0),$I_5$(0),$I_1$(2),$I_5$(2),…$\}$.
Figure 2(a) is similar to Table 1(c), but but \textbf{the orders of a and q are re-arranged}, where
$E^0 \subset \{$6t+1$\}\cup\{$6t+5$\}$ and 
$D^0 \subset \{I_1$(q)$\}\cup\{I_5$(q)$\}$. 
By the re-arranged table in Figure 2(a), SYR($J_0$)=SYR(1)=$\{$1$\} \rightarrow$1, resulting that SYR($J_1$) 
$\rightarrow$1; if SYR($J_1$)$\rightarrow$1, then SYR($J_2$)$\rightarrow$1; if SYR($J_r$)$\rightarrow$1, then SYR($J_{r+1}$)$\rightarrow$1; Thus, SYR($J_s$)$\rightarrow$1, and SYR($D_s$)$\rightarrow$1, $\forall$ s$\in$N.

   
\begin{figure}
\centering
\includegraphics[width=1\linewidth]{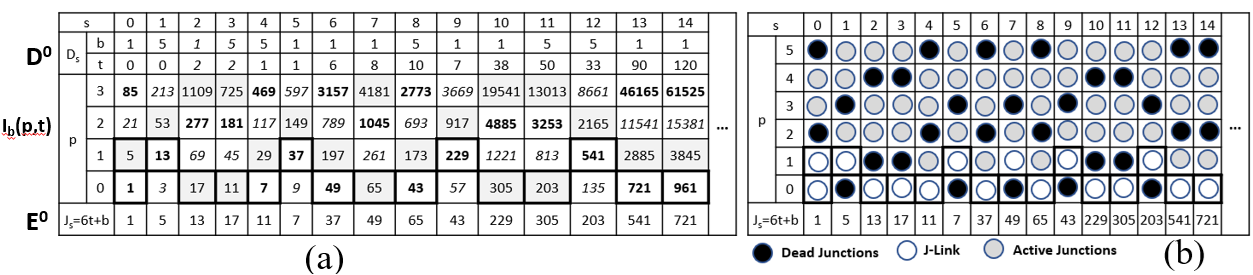}
\caption{\label{fig:IGz}(a) $J_s$ and $D_s$ for $J_0$=1; (b) Black dots, grey dots, and white squares (root candidates).}
\end{figure}

For the Collatz conjecture, as shown in Figure 1(c),  the columns of $\{I_a$(p,q)$\}$, a=1,5, and $\{E_b$(p,q)$\}$, b=1,3,5, are respectively defined as follows,
\begin {equation}
\begin {split}
& I_a(q)=\{m_p\in 2N+1 | Col(m_p)=e_r\}; and  \\
& E_a(q)=\{e_r\in 2N^+ | Col^{r+1}(e_r)=6q+a\}
\end {split}
\end {equation}
where 
\begin {equation}
\begin {split}
& I_a(q)=\{I_a(p,q)\}=\{m_p\},  p,q \in N; and  \\
& E_a(q)=\{E_a(r,q)\}.
\end {split}
\end {equation}
and r=2p+1, if a=1, and r=2p, if a=5.

Given $m_p$=$I_a$(p,q)$\in$2N+1, Col($m_p$)=3*$m_p$+1=$e_r$=2$^{r+1}$x, 
and Col$^{r+1}$($e_r$)= x= 6q+a $\in$2N+1. 
$I_1$(0)=$\{$1,5,21,85,…$\}$=$\{$(4$^{p+1}$-1)/3$\}$, 
$E_1$(0)=$\{$2,4,8,16,…$\}$=$\{2^{r+1}\}$, 
$I_1$(1,0)=5, a=1, p=1, r=2p+1=3, Col($I_1$(1,0))=16=$E_1$(3,0), and 
$E_1$(0,0)=2=2x, x=6q+a, q=0 and a=1. 

For any n$\in$2N+1, SYR(n)=$\{$n=$n_0$,$n_1$,$n_2$,…,$n_s\}$. 
The Collatz sequence COL(n), n$\in$N+, is expressed as 
COL(n)=$R_0 \| R_1 \| R_2 \| \cdots \| R_d \| R_{d+1}$, where
\begin {equation}
\begin {split}
& R_i=\{n_i,E_a(r_i,q_i), E_a(r_{i-1},q_i), …, E_a(0,q_i)\}, and \\	
& n_{i+1}=Col(E_a(0,q_i))=6q_i+a \in 2N+1. 
\end {split}
\end {equation}
Note that $X\|Y$ means the string concatenation of X and Y. 

Let n=$n_0$, if a=1, $r_i$=2$p_i$+1; if a=5, $r_i$=2$p_i$; i=0,1,2,…,d.

•  n=$n_0$=$E_a$($r_0$,$q_0$)=3*$I_{a0}$($p_0$,$q_0$)+1, 
Col$^{r0+1}$($E_a$($r_0$,$q_0$))=6$q_0$+$a_0$=$n_1$=$I_{a1}$($p_1$,$q_1$);

•  $n_1$=$E_a$($r_1$,$q_1$)=3*$I_{a1}$($p_1$,$q_1$)+1, 
Col$^{r1+1}$($E_a$($r_1$,$q_1$))=6$q_1$+$a_1$=$n_2$=$I_{a2}$($p_2$,$q_2$);

• ...

•  $n_d$=$E_a$($r_d$,$q_d$)=3*$I_{ad}$($p_d$,$q_d$)+1, 
Col$^{rd+1}$($E_a$($r_d$,$q_d$))=6$q_d$+$a_d$=$n_{d+1}$=$n_s$=1.

\noindent The sequence COL(n) is
\begin {equation}
\begin {split}
COL(n)=
& \{n=n_0=Col(n_0), Col^2(n_0), …, Col^{r0}(n_0), \\
& n_1, Col(n_1), Col^2(n_1), …, Col^{r1}(n_1), …, \\
& n_d, Col(n_d), Col^2(n_d), …, Col^{rd}(n_d), n_{d+1}\}
\end {split}
\end {equation}
where $n_{d+1}$=$n_s$=1, the sequence COL(n)$\rightarrow$1. 

\section {Tree Trunk for the Syracuse Conjecture}

This section first takes the inverse Syracuse function G in (1.3) to develop an efficient algorithm \textbf{InvSYR} that generates the tree trunk $E^0$=$\{J_0$,$J_1$,$J_2$, …$\}$, main branches $D^0$=$\{D_0$, $D_1$,$D_2$,$D_3$,…$\}$,  the sub-tree trunks, and the sub-branches.  Followed by presenting the properties of the G-functions, where G($J_s$)=$J_{s+1}$, SYR($J_s$)$\rightarrow$1, and SYR($D_s$)$\rightarrow$1, $\forall$s$\in$N.
Note that the inverse Syracuse function and (2.1) assign the positive odd integers to the junctions in $E^0$ and branches $D^0$.

As mentioned in Remark 1.1(4), for any $J_s \in$2N+1, $J_s$=6t+b, $D_{s-1}$=$I_b$(t), this section also develops the algorithm \textbf{PathFinding}, in which the path from n$\in E^r$, down to $E^0$ to verify that SYR(n)$\rightarrow$1, $\forall$n$\in$2N+1.

\subsection {Trunk and Main Branches -- InvSYR Algorithm}

By (1.3),  n=(2$^r$m-1)/3=$I_b$(k,t). For k=0, r=2 if b=1, and r=1 if b=5. Let g=$I_b$(0,t) (mod 6). If $I_b$(0,t)$\neq$3 (mod 6), then G(m)=$I_b$(0,t), otherwise, G(m)=$I_b$(1,t)=4*$I_b$(0,t)+1.

\quad Algorithm \textbf{InvSYR}: s=0, $n_1$=$J_0$=1,
 
\quad \quad \quad Step 1: If $J_0$=1, $D_0$=$I_1$(0), $J_1$=5,  $D_1$=$I_5$(0), s=1;

\quad \quad \quad Step 2: b=mod ($J_s$,6); If b=1, r=2, else (b=5) r=1; m=$2^r J_s$, $m_1$=(m-1)/3;

\quad \quad \quad Step 3: $J_s$=$m_1$; If b=3, then $J_s$=4*$J_s$+1; 

\quad \quad \quad Step 4: s=s+1, b=mod ($m_1$,6), t=($m_1$-b)/6, $D_s$=$I_b$(t), GO to Step 2.

\noindent Example, 

\quad \quad • If $J_0$=1, $D_0$=$I_1$(0), $J_1$=5, $D_1$=$I_5$(0), s=1,

\quad \quad • b=5, r=1, m=10, $m_1$=3, s=2, $J_2$=13=6t+b, b=1, t=2, $D_2$=$I_1$(2), 

\quad \quad • b=1, r=2, m=52, $m_1$=17, s=3, $J_3$=17=6t+b, b=5, t=2, $D_3$=$I_5$(2),

\quad \quad • b=5, r=1, m=34, $m_1$=11, s=4, $J_4$=11=6t+b, b=5, t=1, $D_4$=$I_5$(1),

\quad \quad • …

\noindent where $E^0$=$\{$1,5,13,17,11,…$\}$, and 
$D^0$=$\{I_1$(0),$I_5$(0),$I_1$(2),$I_5$(2),$I_5$(1),…$\}$, as shown in Figure 2(a), 
where, $I_b$(t) can be expressed by Theorem 2.1(2), e.g., b=5, t=2, 
$I_5$(p,t)=[(6t+5)* $4^{p+1}$-2]/6=(17*$4^{p+1}$-2)/6, p$\in$N, 
$I_5$(2)=$\{$11,45,181,725, …$\}$, as shown in Table 1(a).
In Figure 2(b), the entries are marked by black and grey nodes for b=3 and b=1,5, respectively, while the white nodes are for $J_s$. 
In $E^0$, $J_{s+1}$=G($J_s$), s$\in$N, and Syr($J_{s+1}$)=Syr(G($J_s$))=$J_s$. 
Thus, SYR($J_s$)=$\{J_s$,…,$J_1$,$J_0\}$, where SYR$^s$($J_s$)=$J_0$=1, and SYR($J_s$)$\rightarrow$1, i.e., SYR(n)$\rightarrow$1, $\forall$n$\in E^0$.

The function G is defined by the rows, $I_b$(p,t), p=0,1, of the table in Figure 2(a). For the $J_s$, if $I_b$(0,t)$\neq$3 (mod 6), then G($J_s$)=$J_{s+1}$, e.g., $J_2$=13, 
$I_1$(0,2)=17$\neq$3 (mod 6), $J_3$=G($J_2$)=17; and $J_4$=11, $I_5$(0,1)=9=3 (mod 6), $I_5$(1,1)=4*9+1=37, $J_5$=G($J_4$)=37.

\begin {thm} s$\in$N,

(1) If SYR($J_s$)$\rightarrow$1, then SYR($J_{s+1}$)$\rightarrow$1; and 

(2) If SYR($D_s$)$\rightarrow$1, then SYR($D_{s+1}$)$\rightarrow$1.
\end {thm}
\begin {proof}
(1) $E^0$=$\{J_0$,$J_1$,…,$J_s$,…$\}$, by (1.3), G($J_s$)=$J_{s+1}$, Syr($J_{s+1}$)= Syr(G($J_s$))=$J_s$, 
 i.e., SYR($J_s$)=$\{J_s$,...,$J_1$,$J_0\}$,
Syr($J_{s+1}$)=$J_s$, and SYR($J_{s+1}$)=$\{J_{s+1}$,$J_s$,$J_{s-1}$,...,$J_1$,$J_0\}$,
If SYR($J_s$) $\rightarrow$1, then $J_0$=1, and thus SYR($J_{s+1}$)$\rightarrow$1;

\noindent (2) $J_s$=1=6t+b, $D_{s-1}$=$I_b$(t), by (2.2) and (2.3), if SYR($J_s$)$\rightarrow$1, 
then SYR($D_{s-1}$)$\rightarrow$1. 
Let $J_{s+1}$=6t+b, $D_s$=$I_b$(t), if SYR($D_s$)$\rightarrow$1, by (2,2) and (2.3), 
SYR($D_{s+1}$)$\rightarrow$1, by (1), SYR($J_{s+2}$)$\rightarrow$1, and 
SYR($D_{s+1}$)$\rightarrow$1.
\end {proof}

\begin {thm} SYR($J_s$)$\rightarrow$1, and SYR($D_s$)$\rightarrow$1, $\forall$ s$\in$N;
\end {thm}
\begin {proof}
If $J_0$=1, SYR($J_0$)$\rightarrow$1, by Theorem 3.1(1), SYR($J_1$)$\rightarrow$1; 
if SYR($J_1$) $\rightarrow$1, then SYR($J_2$) $\rightarrow$1; 
if SYR($J_s$)$\rightarrow$1, then SYR($J_{s+1}$)$\rightarrow$1, s$\in$N.
Similarly, if SYR($J_1$) $\rightarrow$1, then SYR($D_0$)$\rightarrow$1;
by Theorem 3.1, SYR($D_0$)$\rightarrow$1; 
and if SYR($D_s$)$\rightarrow$1, then SYR($D_{s+1}$)$\rightarrow$1, s$\in$N.
\end {proof}

Let $B_0$ and $F_0$ denote the collection of the active and dead junctions of $D^0$, and $B_{0s}$ and $F_{0s}$ be the collections of the active and dead junctions of $D_s$, respectively. The following theorem derives the expresses of $B_{0s}$ and $F_{0s}$. By Remark 2.1(2), 
let $u^p$=$I_b$(p,t), $v^p$=$u^p$ (mod 6), p=0,1,2. 
If $v^0$=1, then $v^1$=5, and $v^2$=3; if $v^0$=3, then $v^1$=1, and $v^2$=5; and if $v^0$=5, then $v^1$=3, and $v^2$=1; 

\begin {thm} Let $u_p$=$I_b$(p,t), $v_p$=$u_p$ (mod 6), p=0,1,2, 

(1) If $v_0$=1, $B_{0s}$=$\{I_1$(3k+3,t), $I_5$(3k+1,t)$\}$, and $F_{0s}$={$I_3$(3k+2,t)}; 

(2) If $v_0$=3, $B_{0s}$=$\{I_1$(3k+4,t), $I_5$(3k+2,t)$\}$, and $F_{0s}$={$I_3$(3k,t)};  

(3) If $v_0$=5, $B_{0s}$=$\{I_1$(3k+2,t), $I_5$(3k+3,t)$\}$, and $F_{0s}$={$I_3$(3k+1,t)}. 
\end {thm}
\begin {proof}
(1) If $v_0$=1, then $v_1$=5, and $v_2$=3; $F_{0s}$={$I_3$(3k+2,t)}, and 
$B_{0s}$=$\{I_5$(1,t), $I_1$(3k+3,t), 

$I_5$(3k+4,t)$\}$=$\{I_1$(3k+3,t), $I_5$(3k+1,t)$\}$;

\noindent (2) If $v_0$=3, then $v_1$=1, and $v_2$=5; $F_{0s}$={$I_3$(3k,t)}, and $B_{0s}$=$\{I_5$(2,t), $I_1$(3k+4,t), $I_5$(3k+5,t)$\}$  

=$\{I_1$(3k+4,t), $I_5$(3k+2,t)$\}$;

\noindent (3) If $v_0$=5, then $v_1$=3, $v_2$=1; $F_{0s}$={$I_3$(3k+1,t)},
and $B_{0s}$=$\{I_5$(2,t), $I_5$(3k+3,t), $I_1$(3k+5,t)$\}$

=$\{I_1$(3k+2,t), $I_5$(3k+3,t)$\}$
\end {proof}

For s=6, $u_0$=49, $v_0$=1, by Theorem 3.3, $B_{06}$=$\{I_1$(3k+3,t), $I_5$(3k+1,t)$\}$=$\{I_5$(1,6), $I_1$(3,6), $I_5$(4,6), $I_1$(6,6), …$\}$, and
$F_{06}$={$I_3$(3k+2,t)}=$\{I_3$(2,t),$I_3$(5,t),$I_3$(8,t),…$\}$, as shown in Figure 2(b); and 
for s=3, $u_0$=11, $v_0$=5, $B_{03}$=$\{I_1$(2,t), $I_5$(3k+3,t), $I_1$(3k+5,t)$\}$=$\{I_1$(2,3), $I_5$(3,3), $I_1$(5,3), $I_5$(6,3),…$\}$, and 
$F_{03}$={$I_3$(1,t),$I_3$(4,t),$I_3$(7,t),…$\}$.

The algorithm \textbf{InvSYR} generates the tree trunk $E^0$ and the main branches $D^0$,  where $D^0\subset \{I_1$(p,q)$\}\cup\{I_5$(p,q)$\}$. 
SYR($J_s$)$\rightarrow$1 and SYR($D_s$)$\rightarrow$1, $\forall J_s \in E^0$ and 
$D_s \in D^0$; All active junctions in $D^0$ produce the sub-tree trunk $E^1$ and sub-branches $D^1$, where $D^1\subset \{I_1$(p,q)$\}\cup\{I_5$(p,q)$\}$.
SYR($J_s$)$\rightarrow$1 and SYR($D_s$)$\rightarrow$1, $\forall J_s \in E^1$ and 
$D_s \in D^1$;mThe procedure is repeatedly applied, 
for all active junctions in $D^r$, to produce the new trunks $E^{r+1}$ and new sub-branches $D^{r+1}$, where $D^{r+1}\subset \{I_1$(p,q)$\}\cup\{I_5$(p,q)$\}$, and
SYR($J_s$)$\rightarrow$1 and SYR($D_s$)$\rightarrow$1, $\forall J_s \in E^{r+1}$ and $D_s \in D^{r+1}$, , as shown in Figure 1(e). If we assume that the tree trunk, sub-tree trunks, and new sub-branches can grow indefinitely, by (1.6), D=$\cup_{r=0}^{\infty}D_r$=$\{I_1$(p,q)$\}\cup\{I_5$(p,q)$\}$=2N+1.

\subsection {Pathfinding Algorithm}

Given any root junction $J_{a0}$, the algorithm \textbf{InvSYR} generates the junctions $J_{as}$, s$\in$N. 
By (1.7), SYR(n)$\rightarrow$1, $\forall$ n$\in$D=2N+1,
meaning that all positive odd integers always eventually reach the number of 1. 
For example, let $J_{a0}$=1709, the algorithm \textbf{InvSYR} produces the sub-tree trunk $\{$1709,1139,3037,4049,2699,1799,1199,799,4261,5681,…$\}$. SYR(5681) = $\{$5681,4261,799,1199,1799,2699,4049,3037,1139,1709,…$\}$, where n= 5681=$J_{as} \in E^r$, r is unknown for the virtual structure. 
SYR$^9$(5681)=1709=$J_{a0}$. m=(1709-1)/4=427$\neq$3 (mod 6). 
d=(1709-1)/8=53.25$\notin$N, e=(1709-3)/4=106$\in$N, $J_{bs-1}$ =6*106+5= 641$\in E^{r-1}$. 
(By Theorem 2.1(4) and Figure 1(c), Syr(x)=$J_{bs-1}$, $\forall$  x$\in D_{bs-1}$.)
Note that if m=3 (mod 6), m=$J_{bs-1}$ is a dead junction, m$\in E^{r-1}$.  Evidently, $J_{a0}$ is an \textit{invalid root} because that $J_{a0} \notin E^{r-1}$ and $J_{a0}\notin D_{bs-1}$.

By the Syracuse function in (1.2), one can get the sequence Syr$^{26}$(5681)=1, but the Syracuse problem is yet still unsolved. This section develops the algorithm \textbf{PathFinding} to find a path from n to the root junction $J_0$=1 by the virtual structure in Figure 1(f).

First, given n$\in$2N+1, the locations of n, $E^r$, $J_s$, and $D_{s-1}$ are identified: Given n, the parameters b, p, and t are obtained [3], $J_s$=6t+b, $D_{s-1}$=$I_b$(t), 
$J_{s+1}$=$I_b$(p,t), p=0 or 1. For example, n=5681, d=(5681-1)/8=710$\in$N, e=(5681-3)/4=1419.5$\notin$N; $J_s$=5681$\in E^r$, $J_{s-1}$=6*710+1=4261, as shown in Figure 3(c).

\noindent  \textbf{Algorithm PathFinding}: 

Given $J_{as} \in E^r$, to find $J_0$. Let d=($J_{as}$-1)/8, e=($J_{as}$-3)/4, and m=($J_{as}$-1)/4.

Step 0: u=0,

Step 1: If d$\in$N and e$\notin$N, then b=1, t=d, and $J_{as-u}$=6t+b;  

Step 2: If d$\notin$N and e$\in$N, then b=5, t=e, and $J_{as-u}$=6t+b;	

Step 3: If $\notin$N, e$\notin$N, then

\quad \quad  \quad  \quad	if m=3 (mod 6), then $J_{as}$=m, u=u+1, Go to Step 1; 
	 
\quad \quad   \quad   \quad    Else, $J_{as-u}$=6t+b is a valid root, u=0.
        
(If Jas-u is a valid root, $J_{as-u}$=$J_{a0}$ as=u, m$\in E^{r-1}$, and $J_{a0}\in D_{bs-1} \in D^{r-1}$)

\noindent \textbf{Example}: Given $J_{as}$=5681$\in E^r$, to find $J_{a0}$. 

\quad \quad • $J_{as}$=m=5681, d=710$\in$N, e=1419.5$\notin$N, b=1, t=d=710, $J_{as-1}$=6t+b=4261;

\quad \quad • $J_{as-1}$=m=4261, d=532.5$\notin$N, e=1064.5$\notin$N, m=1065=3 (mod 6), 

\quad \quad \quad d=(1065-1)=133$\in$N, e=(1065-3)/4=265.5$\notin$N, b=1, t=d, $J_{as-2}$=799;

\quad \quad •$J_{as-2}$=m=799, d=99.75$\notin$N, e=199$\in$N, b=5, t=d=199, $J_{as-3}$=1199;

\quad \quad •$J_{as-3}$=m=1199, d=149.75$\notin$N, e=299$\in$N, b=5, t=e=299, $J_{as-4}$=1799;

\quad \quad •$J_{as-4}$=m=1799, d=224.75$\notin$N, e=449$\in$N, b=5, t=e=449, $J_{as-5}$=2699;

\quad \quad •$J_{as-5}$=m=2699, d=337.25$\notin$N, e=674$\in$N, b=5, t=e=674, $J_{as-6}$=4049;

\quad \quad •$J_{as-6}$=m=4049, d=506$\in$N, e=1011.5$\notin$N, b=1, t=d=506, $J_{as-7}$=3037;

\quad \quad •$J_{as-7}$=m=3037, d=379.5$\notin$N, e=758.5$\notin$N, m=759=3 (mod 6),

\quad \quad \quad d=94.75$\notin$N, e=189$\in$N, b=5, t=e=189, $J_{as-8}$=1139;

\quad \quad •$J_{as-8}$=m=1139, d=142.5$\notin$N, e=284$\in$N, b=5, t=e=284, $J_{as-9}$=1709;

\quad \quad •$J_{as-9}$=m=1709, d=213.5$\notin$N, e=426.5$\notin$N, m=427$\neq$3 (mod 6), valid root; 

\quad \quad \quad (s-9=0, s=9, $J_{as}$=5681, $J_{a0}$=1709=$I_5$(106)$\in D^r$.)

\quad \quad •$J_s$=m=427$\in$Er-1, d=53.25$\notin$N e=106$\in$N, b=5, t=e=106, $J_{s-1}$=641;

\quad \quad •$J_{s-1}$=m=641, d=80$\in$N, e=159.5$\notin$N, b=1, t=d=80, $J_{s-2}$=481;

\quad \quad •…

   
\begin{figure}
\centering
\includegraphics[width=0.9\linewidth]{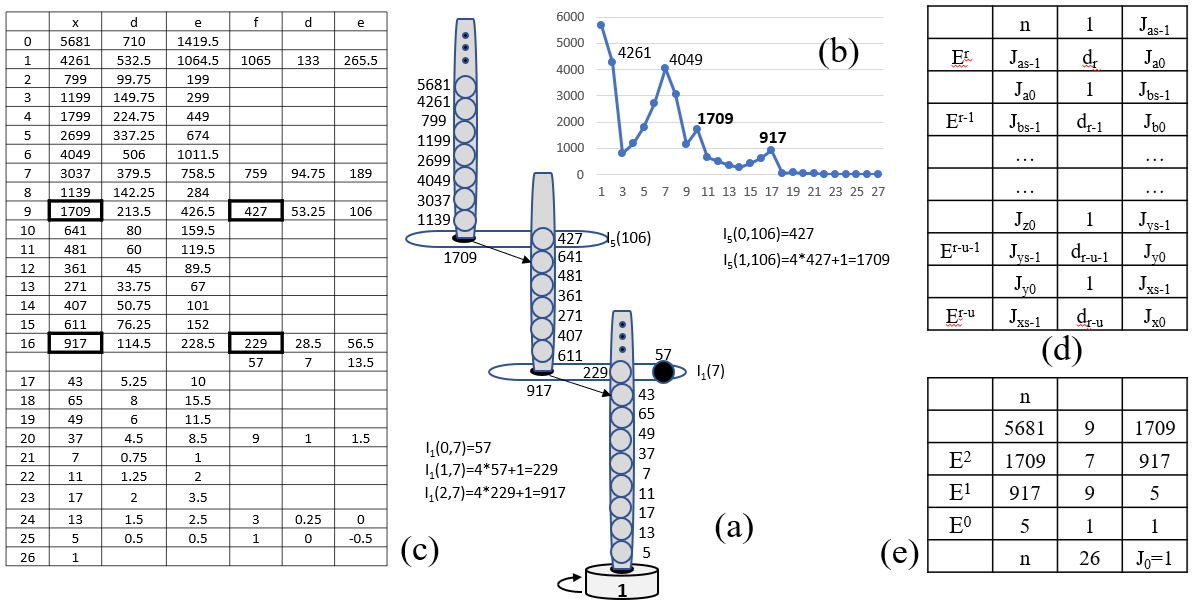}
\caption{\label{fig:IGz}Virtual Tree Trunk for Syracuse Sequence: (a) Tree trunk structure; (b) Plot of SYR(5681); (c) Procedure of searching for the valid roots; (d) Locations of various tree junctions; and (e) Tree trunk and sub-trunks of SYR(5681). }
\end{figure}

For the Syracuse sequence SYR(n)=$\{$n,Syr(n),Syr$^2$(n),…, Syr$^s$(n)$\}$, n$\in$2N+1. 
SYR(n) $\rightarrow$1 if Syr$^s$(n)=1. By (1.3), n=$J_s \in$2N+1, $J_{r-1}$=Syr($J_r$), 
if n=$J_s \in E^0$, then there exists s such that Syr$^s$(n)=Syr$^s$($J_s$)=$J_0$=1. 
However, if $J_{as} \in E^r$, then Syr$^{dr}$(n)=Syr$^{dr}$($J_{as}$)=$J_{a0}$, 
$J_{a0}$ locates at the sub-branch $D_{bs-1}$, as shown in Figure 1(f),  
by Theorem 1.2(4), Syr($J_{a0}$)= $J_{bs-1}$, 
then Syr$^{d(r-1)}$($J_{bs-1}$)=$J_{b0}$. 
The procedure is repeatedly applied, let $J_{z0} \in E^2$,
and locates at the sub-branch $D_{ys-1}$,
Syr$^{d1}$($J_{ys-1}$)=$J_{y0}\in E^1$, 
and at the sub-branch $D_{xs-1}$, and $J_{xs}\in E^0$,
and SYR$^{r0}$($J_{s-1}$)=$J_0$=1. Thus, the path, as shown in Figure 1(f), is
\begin {equation} 
\begin {split}
& n\rightarrow J_{as-1}\rightarrow J_{a0}\rightarrow J_{bs-1}\rightarrow J_{b0}\rightarrow … \\
& \rightarrow J_{z0}\rightarrow J_{ys-1} \rightarrow J_{y0}\rightarrow J_{xs-1}\rightarrow J_{x0}=1
\end {split}
\end {equation}	
as shown in Figure 3(a), and the Syracuse sequence is 
\begin {equation} 
\begin {split}
& SYR(n)=\{n,J_{as-1},J_{as-2},…,J_{a0},J_{bs-1},J_{bs-2},… J_{b0}, ...,J_{z0},  \\
& J_{ys-1},J_{ys-2},… J_{y0},J_{xs-1},J_{xs-2},… J_{x0}=1\}
\end {split}
\end {equation}		
where Syr$^d$(b)=1, where d=d$^r$+d$^{r-1}$+…+d$^1$+d$^0$.
By Figure 3(a), 1709, 917, and 5 are three junctions of $E^2$, $E^1$, and $E^0$, respectively, as shown in Figure 4(e), where Syr$^9$(5681)= 1709, Syr$^7$(1709)=917, Syr$^9$(917)=5, and SYR(5)=1, or Syr$^{26}$(5681)=1, and SYR(5681)$\rightarrow$1.	

Consider the sequence SYR(517), Figure 4(a) shows the tree trunk for demonstrating the sequence SYR(517). Figures 4(c) and 4(d) plot the sequences SYR(517) and COL(517), respectively. Similar to the discussion of SYR(5681) in Figure 3, Figure 4(b) shows that the valid roots are $\{$125,2429,3077,53,5$\}$ which locate at $E^0$, $E^1$, $E^2$, $E^3$, and $E^4$, respectively. Thus, 517=$J_{x5}\in E^4$, and 517=$I_1$(0,86). The path for SYR(517)=$\{$517,97,73,…, 53,5,1$\}$ is listed in Figure 4(e), and SYR(517)$\rightarrow$1. Both 445 and 325 are not located at the sub-branches and they are not the valid roots.

   
\begin{figure}
\centering
\includegraphics[width=0.95\linewidth]{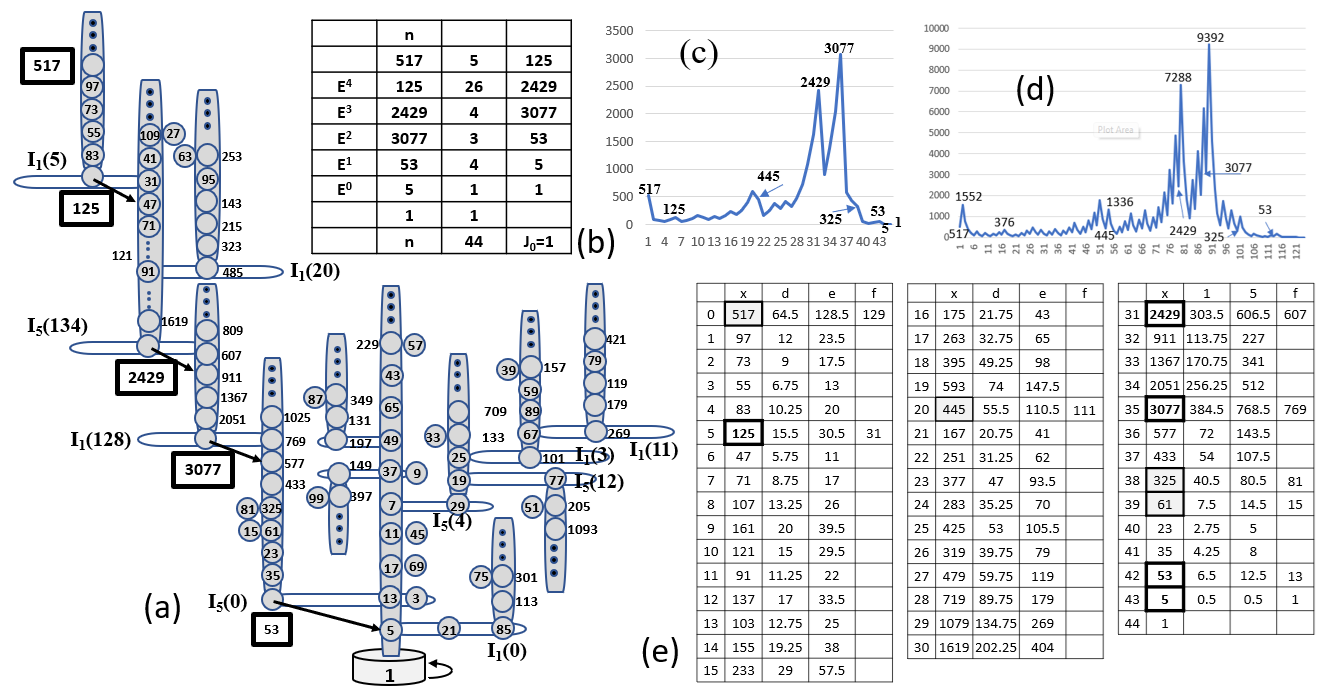}
\caption{\label{fig:IGz}Figure 4: (a) Tree Trunk structure for SYR(517) (b) Tree trunk and sub-trunks of SYR(5681); (c) plot of SYR(517) plot;  (d) plot of Collatz sequence COL(517), and (e) Procedure of searching for the valid roots. .}
\end{figure}

If $J_{x0}$ locates at the sub-branch $D_{bs-1}$, then $J_{x0}$ is a valid root.
 On the other hand, all junctions $J_{y0}$ in $E^r$, $J_{ys}\neq$3 (mod 6), $J_{x0}$ is an invalid root means that it does not locate at any sub-branch in $D^{r-1}$.

Figure 4(b) shows the properties of SYR(517), where n=517$\in E^r$, r is unknown, 125$\in D^{r-1}$, and Syr$^5$(517)=125, Syr$^{26}$(125)=2429, Syr$^4$(2429)=3077, Syr$^7$(3077)=53, Syr(53)=5, Syr(5)=1, and Syr$^{44}$(517)=1. Thus, s=44, Syrs(n)=1, SYR(517)$\rightarrow$1. 

Interestingly, when the five valid roots, $\{$125,2429,3077,53,5$\}$, and some other valid roots $\{$29,53,77,85,101,149,197,269,485,2429$\}$, are used, as shown in Figure 3(a), by Algorithm \textbf{InvSYR},  the tree-trunk and sub-tree trunks include the 50 odd positive integers ranged from 1 to 99. 
By Theorem 3.2, SYR(n)$\rightarrow$1, n=2k+1, k=0$\sim$49. Based on the simulation results, the roots produce a total of 155 positive odd integers which are less than 500, i.e., 155/250=62\%, among them, it has 50, 41, 26, 19, and 19 are within the ranges, 1$\sim$99, 101$\sim$199, 201$\sim$299, 301$\sim$399, and 401$\sim$499, respectively. 

Further, if the valid roots $\{$173,221,245,293,317,341,413$\}$, are added, the total numbers of the odd integers less than 500 become 176, 
i.e., 176/250=70.4\%. Among them, the numbers of the odd integers are 50, 50, 34, 22, and 20 in the same ranges, respectively. This results that SYR(n)$\rightarrow$1, for n=2k+1, k=0$\sim$99. As the number of valid roots are further applied, the value k will significantly increase. By (1.5), 
D=$\cup_{r=0}^{\infty}D_r$=2N+1,
meaning that all positive odd integers always eventually reach the number of 1.

\subsection {Proof of the  Syracuse Conjecture}

By (1.6), D=$\cup_{r=0}^{\infty}D_r$=$\{I_1$(p,q)$\}\cup\{I_5$(p,q)$\}$=2N+1.
By (1.7), SYR(n)$\rightarrow$1, $\forall$n$\in$2N+1, proves the Syracuse conjecture.
Let $J_{x0}$ and $J_{y0}$ be two arbitrary roots of $E^a$ and $E^b$, respectively, where a$\neq$b. If $J_{x0}\neq J_{y0}$, then $\{J_x\}\cap\{J_y\}$=$\phi$ and $D^a\cap D^b$=$\phi$,  
meaning that all junctions of $E^a$ and $E^b$ are distinct, and all junctions of $D^a$ and $D^b$ are also distinct.

\begin {lemma} Let $J_0$=1 produce $E^0$ and $D^0$, and $J_a$, $J_b\in E^0$, $D_a$, $D_b\in D^0$, if a$\neq$b, then $J_a \neq J_b$ and $D^a\cap D^b$=$\phi$.
\end{lemma}
\begin {proof}
Suppose that $D_x\cap D_y \neq \phi$, i.e., both columns in $\{I_a$(p,q)$\}$, a=1,5, of Table 1(a), are not distinct which contradicts to Theorem 1.1(3). 
Thus, $D^a\cap D^b$=$\phi$. 
If $J_a$=$J_b$, then $J_{a+1}$=$J_{b+1}$, and $D^a$=$D^b$, this contradicts to $D^a\cap D^b$=$\phi$. 
\end {proof}

For any root $J_{x0}\in E^r$, by Lemma 3.4,  $J_{xa} \neq J_{xb}$ and $D^a\cap D^b$=$\phi$. If $J_{x0}\neq J_0$, then $\{J_x\}\cap\{J_y\}$=$\phi$, and $D_x \cap D_y$=$\phi$,  as shown in the following theorems.

\begin {lemma} If $J_{x0}\neq J_0$, then $\{J_x\}\cap\ E^0$=$\phi$, and $D_x \cap D^0$=$\phi$ 
\end{lemma}
\begin {proof}
Suppose that $\{J_x\}\cap\ E^0 \neq \phi$, there exists $J_z \in \{J_x\}\cap\ E^0$,
 by Algorithm \textbf{RJDGen}, given $J_z$, both roots $J_{x0}$ and $J_0$ are generated, respectively, and thus $J_{x0}$=$J_0$ 
which contradicts to $J_{x0} \neq J_0$. Thus, $\{J_x\}\cap\ E^0$=$\phi$.
By Lemma 3.4, if $I_b$(p,t)$\in D_x$ and $I_b$(p,t)$\in D^0$,
the roots are $J_{x0}$ and $J_0$, respectively, and $J_{x0}$=$J_0$, 
which contradicts to $J_{x0} \neq J_0$=1. Thus, $D_x \cap D^0$=$\phi$, 
\end {proof}

\begin {lemma} If $J_{x0}\neq J_{y0}$, then  $\{J_x\}\cap\{J_y\}$=$\phi$ and $D_x \cap D_y$=$\phi$ 
\end{lemma}
\begin {proof}
The proof is similar to that of Lemma 3.5. 
\end {proof}

By Theorem 3.1, if SYR($J_s$)$\rightarrow$1, then SYR($J_{s+1}$)$\rightarrow$1. Thus, If SYR($J_{s+1}$)$\nrightarrow$1, then SYR($J_s$)$\nrightarrow$1; and if SYR($D_{s+1}$)$\nrightarrow$1, then SYR($D_s$)$\nrightarrow$1. The following theorem proves (1.5).

\begin {thm}
D=$\cup_{r=0}^{\infty}D^r$=$\{I_1$(p,q)$\}\cup\{I_5$(p,q)$\}$=2N+1.
\end {thm}
\begin {proof} (1) \textbf{$\{I_1$(p,q)$\}\cup\{I_5$(p,q)$\}\subseteq$D=$\cup_{r=0}^{\infty}D^r$}

\noindent  For any $D_x$=$I_b$(t)$\in \{I_1$(p,q)$\}\cup\{I_5$(p,q)$\}$,
$J_x$=6t+b, suppose $D_x \notin$D=$\cup_{r=0}^{\infty}D_r$, 
without loss of generality, let $D_x \notin 
\cup_{i=r}^{\infty}D^i$,  $D_x \notin D^r$ implies that $J_x$=6t+b$\notin E^r$,
or SYR($J_x$)$\nrightarrow$1.
(If $J_x$=6t+b$\in E^r$ and SYR($J_x$)$\rightarrow$1, then $J_x \in E^r$ and 
$D_x \in D^r$.)
If $J_x$=6t+b$\notin E^r$, SYR($J_x$)$\nrightarrow$1,  
 then SYR($J_{x-1}$)$\nrightarrow$1, and SYR($J_{x0}$)$\nrightarrow$1, as shown in Figure 1(f), Syr($J_{x0}$)= $J_{x0-1} \in E^{r-1}$ and SYR($J_{xs}$)$\nrightarrow$1, 
$\forall J_{xs} \in E^{r-1}$. The procedure is repeatedly applied until $J_{z0}$ is obtained and located at $D_{ys-1}$. 
Syr($J_{z0}$ )=$J_{ys-1} \in E^1$, and SYR($J_{ys}$)$\nrightarrow$1,
$\forall J_{ys} \in E^1$. $J_{y0}$ locates at $D^s$, 
Syr($J_{y0}$ )=$J_{ys-1} \in E^0$, and SYR($J_s$)$\nrightarrow$1,
$\forall J_s \in E^0$, and $J_0$=1 is terminated, where 
SYR($J_0$)=SYR(1)$\nrightarrow$1, contradicting that SYR($J_0$)=SYR(1)$\rightarrow$1. Therefore, $D_x \in$D, and $\{I_1$(p,q)$\}\cup\{I_5$(p,q)$\}\subseteq$D

\noindent (2) \textbf{D$\subseteq \{I_1$(p,q)$\}\cup\{I_5$(p,q)$\}$}

\noindent  For any $D_x$=$I_b$(t)$\in \cup_{r=0}^{\infty}D^r$, 
$I_b$(t)$in \{I_1$(p,q)$\}\cup\{I_5$(p,q)$\}$, i.e., 
D$\subseteq \{I_1$(p,q)$\}\cup\{I_5$(p,q)$\}$.

\noindent  Thus, by (1) and (2), $\{I_1$(p,q)$\}\cup\{I_5$(p,q)$\}$=D=$\cup_{r=0}^{\infty}D^r$
\end {proof}

\begin {thm} (Syracuse Conjecture) SYR(n)$\rightarrow$1, $\forall$n$\in$2N+1.
\end {thm}
\begin {proof}
By Theorem 3.2, SYR($J_s$)$\rightarrow$1, $\forall J_s \in E^0$, and 
SYR($D_s$)$\rightarrow$1, $\forall J_{x0} \in B_0$.
Similarly, SYR($J_s$)$\rightarrow$1, $\forall J_s \in E^1$, and 
SYR($D_s$)$\rightarrow$1, $\forall J_{y0} \in B_1$.
SYR($J_s$)$\rightarrow$1, $\forall J_s \in E^r$, and 
SYR($D_s$)$\rightarrow$1, $\forall J_{z0} \in B_r$.
and so on. This concludes that SYR($J_x$)$\rightarrow$1, $\forall J_x \in$E=$\cup_{r=0}^{\infty}E^r$, and by Theorem 3.7, SYR($D_x$)$\rightarrow$1,
for all root junctions in
D=$\cup_{r=0}^{\infty}D^r$=$\{I_1$(p,q)$\}\cup\{I_5$(p,q)$\}$ =2N+1.
meaning that, n=$D_s$=$I_b$(t) and SYR(n)$\rightarrow$1, 
$\forall$n$\in \{I_1$(p,q)$\}\cup\{I_5$(p,q)$\}$=2N+1.
\end {proof}

\section {The Proof of the Collatz Conjecture}

\subsection {Trunk and Main Branches -- \textbf{InvCOL} Algorithm}

By (1.1) (Collatz function), m=Col(n)=3n+1, if n$\in$2N+1, and m=Col(n)=n/2 if n$\in N^+$. Let d=(m-1)/3, if n$\in$2N+1, and d$\neq$3 (mod 6), by (1.4), H(m)=d; otherwise H(m)=2m. 
The H-sequence is V$^0$=$\{h_0$,$h_1$,…,$h_s$,…$\}$, $h_r\neq$3 (mod 6).

\noindent  \textbf{Algorithm InvCOL}:

Step 1: $h_0$=1, $h_1$=2, $h_2$=4 (trivial cycle), s=2,

Step 2: $h_{s+1}$=2*$h_s$, d=($h_{s+1}$-1)/3 (mod 6), s=s+1

Step 3: If d$\in$2N+1, and  d$\neq$3, then  s=s+2, $h_s$=d, Go to Step 2.

\noindent \textbf{Example 1}: $h_{s+1}$=2*$h_s$, d=($h_{s+1}$-1)/3 (mod 6)

•	$h_0$=1, h1=2, $h_2$=4; s=2

•	$h_3$=8, d=7/3$\notin$2N+1; s=s+1=3, 

•	$h_4$=16, d=5, s=s+1=4, d$\notin$2N+1, and d$\neq$3, $h_5$=5, s=s+1=5;

•	$h_6$=10, d=3, s=s+1=6, d=3, 

•	$h_7$=20, d=19/3$\notin$2N+1; s=s+1=7, 

•	$h_8$=40; d=13$\in$2N+1, and d$\neq$3; s=s+2=9, $h_s$=$h_9$=d=13;

•	$h_{10}$=2*13=26; d=25/3$\notin$2N+1; s=10

•	…

\noindent $V^0$=$\{$1,2,4,8,16,5,10,20,40,13,26,….$\}$, as shown in Figure 1(c). $h_0$=1, 
and H($h_s$)=$h_{s+1}$, 

\quad \quad \quad COL(26)=$\{$26,13,40,20,10,5,16,8,4,2,1$\}$.

\noindent  In $V^0$, $h_{s+1}$=H($h_s$), s$\in$N, and Col($h_{s+1}$)=Col(H($h_s$))=$h_s$. Thus, COL($h_s$)=$\{h_s$,$h_{s-1}$,…, $h_1$,$h_0\}$, where $h_0$=26, Col$^{10}$($h_0$)=Col$^{10}$(26)=1, and COL($h_s$)$\rightarrow$1, i.e., COL(n)$\rightarrow$1, $\forall$n$\in V^0$.

\noindent \textbf{Example 2}: $h_{x0}$=125, s=0

•	$h_{x0}$=125, d=124/3$\notin$2N+1;
 
•	$h_{x1}$=250, d=249/3=83, d$\in$2N+1, d$\neq$3, $h_{x2}$=83;

•	$h_{x3}$=166, d=165/3=55, d$\in$2N+1, d$\neq$3, $h_{x4}$=55;

•	$h_{x5}$=110, d=99/3=33=3 (mod 6);  

•	$h_{x6}$=220, d=219/3=73, d$\in$2N+1, d$\neq$3, $h_{x7}$=73;

•	$h_{x8}$=146, d=145/3$\notin$2N+1;

•	$h_{x9}$=292, d=291/3=97, d$\in$2N+1, d$\neq$3, $h_{x10}$=97;  

•	$h_{x11}$=194, d=193/3$\notin$2N+1;

•	$h_{x12}$=388, d=387/3=129=3 (mod 6);
 
•	$h_{x13}$=776, d=575/3$\notin$2N+1

•	$h_{x14}$=1552, d=1551/3=517, d$\in$2N+1, d$\neq$3; $h_{x15}$=517;

•	…

\subsection {Pathfinding Algorithm}

Similar to the \textbf{PathFinding} algorithm discussed in Section 3.2, the virtual tree trunk structure in Figure 1(f), are used to find the path for the Collatz sequences. Figure 4(d) plots the Collatz sequence COL(517).

Let $h_0$=1, and n=517$\in N^+$, by Algorithm  \textbf{InvCOL}, $h_{x0}$=125 generates the tree trunk $V^0$=$\{h_{x0}$,$h_{x1}$,$h_{x2}$,…$\}$, main branches,  the sub-tree trunks, and the sub-branches

\quad \quad $V^0$=$\{$\textbf{125},250,$\textbf{83}$,166,$\textbf{55}$,110,220,$\textbf{73}$,146,292,
$\textbf{97}$,194,388,776,1552,$\textbf{517}$,…$\}\subset E^r$, 

\noindent  and H$^{15}$(125)=517; 

\quad \quad  COL(517)=$\{\textbf{517}$,1552,776,388,194,$\textbf{97}$,292,146,$\textbf{73}$,220,110,$\textbf{55}$,166,$\textbf{83}$,250,\textbf{125}$\}$, 

\noindent   and Col$^{15}$(517)=125. Let
n=$h_{x15}$=517$\in E^r$ is with the valid root $J_{a0}$=125$\in I_1$(5). Col($J_{a0}$) =Col(125)=376=$J_{as-1}$. Similarly, 
$H^{63}$(2429)=125, Col$^{63}$(125)=2429, $J_{b0}$= 2429$\in I_5$(134), and Col($J_{b0}$)=Col(2429)=7288=$J_{as-1} \in E^{r-1}$; 
H$^{10}$(3077)=2429, Col$^{10}$(2429) =3077, $J_{b0}$= 3077$\in I_1$(128),
and Col($J_{c0}$)=Col(3077)=9232=$J_{cs-1} \in E^{r-2}$;
H$^{26}$(53)=3077, Col$^{26}$(3077)= 53, $J_{d0}$=53$\in I_5$(0),
and Col($J_{d0}$)=Col(53)=160=$J_{ds-1} \in E^{r-3}$;
H$^6$(5)=53, Col$^6$(53)=5, $J_{d0}$=5$\in I_1$(0),
and Col($J_{e0}$)=Col(5)=10$J_{es-1} \in E^{r-4}$;
H$^5$(1)=5, Col$^5$(5)=1=$J_0$. $E^{r-4}$=$E^0$, r-4=0, or r=4, n=517$\in E^4$. 
H$^{125}$(1)=517; and Col$^{125}$(517)=1, thus, COL(517)$\rightarrow$1. It should be mentioned that $\{$125,2429,3077,53,5$\}$ are the \textit{valid roots} and locate at $E^4$, $E^3$, $E^2$, $E^1$, and $E^0$, respectively.

The \textbf{InvCOL} algorithm generates $E^0$ and $D^0$ by the root junction $J_0$=1. Figure 5(b) the virtual tree trunk for SYR(169), 
where the procedure of searching for the valid roots is shown in Figure 5(c) and Figure 5(e) plots the sequence SYR(169).  
Based on the tree trunk and subbranches, Figure 5(d) shows the virtual structure, and Figure 5(f) is the plot of COL(169). 
For the sequence COL(45), where n=45=3 (mod 6), n$\notin E^0$ and locates at a dead junction. 45=$I_5$(1,2), by the Collatz function, 
COL(45)=$\{$45,136,68,34,17$\}$, by Figure 5(a), COL(17)=$\{$17,52,26,13,40,20,10,5,16,8,4,2,1$\}$. 
Thus, COL16(45)=1, and COL(45)$\rightarrow$1. The junctions $\{$77,29,1$\}$ are the valid roots and locate at $E^2$, $E^1$, and $E^0$, respectively.

   
\begin{figure}
\centering
\includegraphics[width=0.95\linewidth]{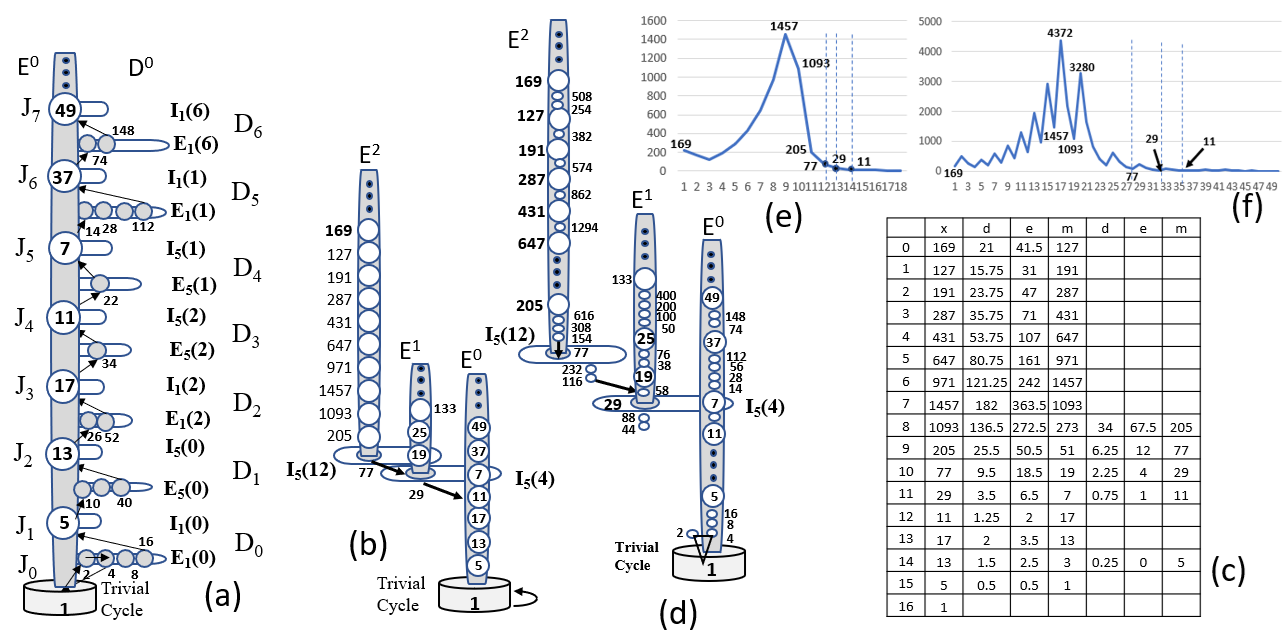}
\caption{\label{fig:IGz}Figure 4: (a) The InvCOL algorithm generates both D$^0$ and E$^0$ by the root junction $J_0$=1; (b) Virtual tree trunk for SYR(169); (c) Procedure of searching for the valid roots; (d) Virtual tree trunk for COL(169); (e) Plot of SYR(169); and (f) Plot of COL (169).}
\end{figure}

\subsection {Proof of the Collatz Conjecture}

By Theorem 3.8, the Syracuse conjecture holds. The proofs of the following two theorems prove the Collatz conjecture.

\begin {thm} If SYR(n)$\rightarrow$1, then COL(n)$\rightarrow$1, 
$\forall$n$\in\{I_1$(p,q)$\}\cup\{I_5$(p,q)$\}$=2N+1.
\end {thm}
\begin {proof}
For any n$\in$2N+1, by (2.7), COL(n)=$\{$n=$n_0$=Col($n_0$), Col$^2$($n_0$), …, Col$^{r0}$($n_0$), $n_1$, Col($n_1$), Col$^2$($n_1$), …, Col$^{r1}$($n_1$), ...,
 $n_d$, Col($n_d$), Col$^2$($n_d$), …, Col$^{rd}$($n_1$), $n_{d+1}\}$. 
and SYR(n)=$\{$n=$n_0$, $n_1$, …, $n_d$, $n_{d+1}\}$. 
If SYR(n)$\rightarrow$1, then $n_{d+1}$=1, and COL(n)$\rightarrow$1. 
\end {proof}

\begin {thm} If COL(n)$\rightarrow$1, $\forall$ n$\in\{I_1$(p,q)$\}\cup\{I_5$(p,q)$\}$, 
then COL(n)$\rightarrow$1, 
$\forall$ n$\in$($\{E_1$(p,q)$\}\cup\{E_3$(p,q)$\}\cup\{E_5$(p,q)$\}$).
\end {thm}
\begin {proof}
If $m_p \in$2N+1=$\{I_1$(p,q)$\}\cup\{I_5$(p,q)$\}$, 
by (2.4) and (2.5), Col($m_p$)=$e_r$=$E_b$(r,q)=$2^{r+1}$x, x=6q+b. 
If COL($m_p$)$\rightarrow$1, then COL($e_r$)$\rightarrow$1, or COL($E_b$(r,q))$\rightarrow$1, b=1,3,5, i.e., COL(n)$\rightarrow$1, 
\noindent $\forall$ n$\in$($\{E_1$(p,q)$\}\cup\{E_3$(p,q)$\}\cup\{E_5$(p,q)$\}$).
\end {proof}

By Theorem 3.8, SYR(n)$\rightarrow$1, $\forall$ n$\in\{I_1$(p,q)$\}\cup\{I_5$(p,q)$\}$=2N+1. Thus, by Theorem 1.1, COL(n)$\rightarrow$1, $\forall$ n$\in\{I_1$(p,q)$\}\cup\{I_5$(p,q)$\}$=2N+1.
By Theorem 1.2, COL(n)$\rightarrow$1, 
$\forall$ n$\in$($\{E_1$(p,q)$\}\cup\{E_3$(p,q)$\}\cup\{E_5$(p,q)$\}$).
This concluded that COL(n)$\rightarrow$1, 
$\forall$ n$\in$($\{I_1$(p,q)$\}\cup\{I_5$(p,q)$\}$)
$\cup$($\{E_1$(p,q)$\}\cup\{E_3$(p,q)$\}\cup\{E_5$(p,q)$\}$)=(2N+1)$\cup$(2N$^+$)=N$^+$. \textit{The Collatz conjecture holds}.

\begin {remark} Proof of the Syracuse and Collatz Conjectures
\begin{enumerate}
\item Both Collatz and Syracuse conjectures can be described conceptually by the tree trunk, as shown in Figure 1. The junctions of the tree trunk and sub-tree trunks are built up by the developed inverse Collatz (Syracuse) functions. Each junction produces a branch. The junctions in all produced tree trunk and sub-tree trunks are distinct, so are the main branches and sub-branches, proved by Lemmas 3.4-3.6. Conceptually assuming that the trunk and branches can grow indefinitely, by (1.6) and Theorem 3.7, the junctions of the branches includes all n$\in$2N+1 for the Syracuse sequences, and all n$\in N^+$ for the Collatz sequences;
\item Based on the virtual structure of Figure 1(f), given n, the Algorithm \textbf{PathFinding} finds a path starting from the given n down to the root junction $J_0$=1. If n=$J_{as}$ locates at $E^r$ or $D_{as-1}$ of $D^r$, the next step is $J_{as-1}$, and down to $J_0$ the root of the sub-tree trunk which locates at a $D_{bs-1}$,
and $J_{bs}$ locates at $E^{r-1}$. The curve of the junctions between $J_{as-1}$ to $J_{a0}$ may be up and down, but, by Figure 1(e), the trunks definitively act likes the stairs starting from $E^r$ down to $J_{xs-1}$ of $E^0$ and then $J_0$=1 to prove the Collatz (Syracuse) conjectures.
\item The Syracuse conjecture is proved by Theorem 3.8, thus, by Theorems 4.1 and 4.2, Collatz conjecture is proved.

\end{enumerate}
\end {remark}


\begin{quote}
\begin{small}
 \noindent \textsc{Chin-Long Wey}.\\
        Department of Electronics and Electrical Engineering \\
        National Yang Ming Chiao Tung University  \\
        Hsinchu, Taiwan  \\
        E-mail: \texttt{wey@nycu.edu.tw}\\
	ORCID: 0009-0002-7063-6189

\end{small}
\end{quote}


\begin{thebibliography}{1}

{
\bibitem{Lagarias}
J.C. Lagarias (2010). \textit{The 3x+1 problem: an Overview}.
arXiv preprint arXiv: 2111.02635 (2010).
}
\bibitem{Tao2019}
T. Tao, (2019). \textit{Almost all orbits of the Collatz map attain almost bounded values}.
arXiv preprint arXiv: 1909.03562 (2019).

\bibitem{Wey}
C.L. Wey (2023). \textit{Proof of Collatz conjecture by Collatz Graph}.
arXiv preprint arXiv: 2309.09991v2 (2024)


\end{thebibliography}
\end{document}